\documentclass{amsart}
\usepackage{hyperref}
\hypersetup{colorlinks=false}

\usepackage{pgf}
\usepackage{tikz}
\usetikzlibrary{cd}

\hypersetup{colorlinks=false}

\newtheorem{theorem}{Theorem}[section]
\newtheorem{prop}{Proposition}[section]
\newtheorem{lemma}{Lemma}[section]
\newtheorem{corollary}{Corollary}[section]

\DeclareMathOperator{\id}{id}
\DeclareMathOperator{\clos}{cl}
\DeclareMathOperator{\interior}{int}
\DeclareMathOperator{\Inv}{Inv}

\newcommand{\w}[1]{\widetilde #1}

\begin{document}

\title[On determining the homological Conley index]{On determining the homological Conley index of
Poincar\'e maps in autonomous systems}

\author{Roman Srzednicki}

\address{
Institute of Mathematics\\ 
Faculty of Mathematics and Computer Science\\
Jagiellonian University\\
ul. {\L}ojasiewicza 6\\
30-348 Krak\'ow, Poland}

\email{srzednicki@im.uj.edu.pl}

\dedicatory{Dedicated to the memory of Professor Andrzej Granas}

\subjclass[2020]{Primary: 37B30; Secondary: 37B35}

\keywords{Poincar\'e map; isolated invariant set; index pair; Conley index}

\thanks{This research is partially supported by the Polish National Science Center
under Grant No. 2014/14/A/ST1/00453}

\begin{abstract}%
A theorem on computation of the homological Conley index of an isolated
invariant set of
the Poincar\'e map associated to a section in a rotating local dynamical system $\phi$ is proved.
Let $(N,L)$ be an index pair for a discretization $\phi^h$ of $\phi$, where $h>0$, and
let $S$ denote the invariant part of $N\setminus L$; it follows that
the section $S_0$ of $S$ is an isolated invariant set of 
the Poincar\'e map. 
The theorem asserts that if 
the sections $N_0$ of $N$ and $L_0$ of $L$ are ANRs,
the homology classes $[u_j]$ of some 
cycles $u_j$ form a basis of $H(N_0,L_0)$,
and for some scalars $a_{ij}$, the cycles $u_j$ and $\sum a_{ij}u_i$ are
homologous in the covering pair $(\widetilde N,\widetilde L)$ of $(N,L)$ and the homology
relation is preserved in $(\widetilde N,\widetilde L)$ under the transformation induced by $\phi^t$ for $t\in [0,h]$
 then the homological Conley index of $S_0$
is equal to the 
Leray reduction of the matrix $[a_{ij}]$. In particular, no information on the values of the Poincar\'e map or its approximations is
required.
In a special case of the system generated by a $T$-periodic non-autonomous ordinary differential
equation with rational $T/h>1$, the theorem was proved in the paper M.\,Mrozek, R.\,Srzednicki, and F.\,Weilandt, 
SIAM J. Appl. Dyn. Syst. 14 (2015), 1348-1386, and it motivated a construction of an algorithm for
determining the index.
\end{abstract}

\maketitle

\section{Introduction}
\label{sec:intro}

Let $v$ be a smooth vector-field on a manifold $X$. Assume that
there is a smooth map $\theta$ from $X$ to the circle, represented as $\mathbb R/\mathbb Z$, 
such that the derivative of $\theta$ at $v(x)$ is positive for every $x$.
It follows that the local dynamical system $\phi$ generated by $v$ is rotating, which roughly means
that for every $x\in X$
the angle function
$\theta(\phi(x,\cdot))$ is increasing. 
The transformation which sends a point 
$x$ of a section 
$X_a:=\theta^{-1}(a)$ to
the point of its return to $X_a$ when the function
$\theta(\phi(x,\cdot))$ makes a full rotation around the circle
is called a Poincar\'e map; in the case $a=0$ we denote it by 
$\Pi$.  The dynamics of $\Pi$ determines  the
qualitative properties of $\phi$; in particular, each periodic 
point of $\Pi$ is an initial point of
a periodic trajectory of $\phi$ and chaotic dynamics of $\Pi$ encodes
the corresponding chaotic behavior of trajectories of $\phi$.
\par
In the present paper we focus our interest on isolated invariant sets
of $\phi$; each such set 
is also an isolated invariant set for the discretization $\phi^t=
\phi(\cdot,t)$ for every $t>0$  and vice versa: an isolated invariant
set $S$ of $\phi^t$ is an isolated invariant set of $\phi$. Moreover,
$S_0=S\cap X_0$ is an isolated invariant set of the map $\Pi$ (cf. Proposition~\ref{prop:s_invariant}). 
\par
Both in the continuous-time setting (like $\phi$ above) and the discrete-time 
setting (i.e. the iterations of a map),
properties of isolated invariant sets are investigated using their 
index pairs.
An index pair $(N,L)$ for
an isolated invariant set $S$ is a pair of compact 
sets such that  $S$ is contained in the interior of $N\setminus L$,
the closure of $N\setminus L$ does not contain any larger invariant set then $S$, $L$ is positively
invariant in $N$, and $L$ is an exit set of $N$. If $f$ is a discrete-time local dynamical system,
 the induced self-map $f_{(N,L)}$ of $N/L$ is continuous. 
The set of non-zero eigenvalues
of the linear map
$H(f_{(N,L)})$, i.e. the map induced by $f_{(N,L)}$ in homologies with coefficients in a field,  
does not depend of the choice of an index pair $(N,L)$. If the domain
of $f$ is an ANR, the non-zero eigenvalues determine 
the Lefschetz numbers of $f_{(N,L)}$ and its iterates, and provide the fixed point index and the Lefschetz
zeta function of $f$ in a neighborhood of $S$, hence they enable 
to detect fixed and periodic points of $f$, and, moreover,
in the case of the smooth manifold, to estimate the entropy of $S$. 
Actually, the mentioned invariance of the eigenvalues follows from the invariance of the homological Conley
index $CH(S,f)$ defined as the
conjugacy class of the Leray
reduction $R$ of $H(f_{(N,L)})$. In Section~\ref{sec:basic} we recall basic notions related to the homological
Conley index; for a complementary information we refer to \cite{baker,m-fm,m-tams,m-etds}.
\par
The above remarks indicate that some important facts on
the dynamics 
generated by the vector-field $v$ in 
a neighborhood of an isolated invariant set $S$ 
can be detected by the index $CH(S_0,\Pi)$. 
However, in practice it is difficult to determine 
the index directly, even though it
is stable with respect to small perturbations of the map $\Pi$. 
An essential obstacle is caused by the exponential growth of errors of numerical 
methods for trajectories of $\phi$; if the first return time to $X_0$ is large, usually it is impossible
to get an approximation of $\phi$ with a required accuracy. There might be other obstacles; for example, 
in systems with a strongly expansive direction the set of points
such that the Poincar\'e map is defined can be too small to be detected
since a vast majority of points do not return 
to the section and escapes to infinity. On the other hand, for a small enough step-size $h$, numerical methods provide 
satisfactory approximations of the discretization $\phi^h$ and there are computation algorithms
which provide index pairs for discrete-time dynamical systems in a rigorous way (see \cite{ft,m-ipa,msw}). 
\par
We show that a suitable information
on an index pair $(N,L)$ with respect to a discretization $\phi^h$ is sufficient to get the value of $CH(S_0,\Pi)$,
without computing $\Pi$ or its approximations.
In Theorem~\ref{thm:main}, being 
the main result of the present paper,  
a sequence $u_j$, $j=1,\ldots,n$, of cycles of the pair $(N_0,L_0)$ is considered, where $N_0=N\cap X_0$ and $L_0=L\cap X_0$. It is assumed that $N_0$ and $L_0$ are ANRs, the homology classes $[u_j]$ 
form a basis of the graded vector space $H(N_0,L_0)$, 
and in the covering pair $(\widetilde N,\widetilde L)$ of $(N,L)$ over $[0,1]$, each $u_j$
over $0$
is homologous to a linear combination $\sum a_{ij}u_i$ over $1$
and the homology relation is preserved  in $(\widetilde N,\widetilde L)$ under 
the transformation induced by $\phi^t$ 
for $t\in [0,h]$. (A precise formulation of the latter assumption: $\left(u_j,\sum a_{ij}u_i\right)$
is an $h$-movable pair of contiguous cycles.)
Theorem~\ref{thm:main} asserts
that $CH(S_0,\Pi)$ is equal to the Leray reduction of the matrix $[a_{ij}]$.
\par
As a consequence of the theorem we prove that under the above assumptions on $(N,L)$
and $(N_0,L_0)$, if a homotopy 
$F\colon (N_0,L_0)\times [0,1]\to (N,L)$ satisfies $F(x,0)=x$ and 
$F(x,a)\in X_a$ for all $x\in N_0$ and $a\in [0,1]$, and  
$\phi^t(F(N_0,a))\subset N$ and $\phi^t(F(L_0,a))\subset L$ for all $a\in [0,1]$ and $t\in [0,h]$
then $CH(S_0,\Pi)$ is equal to the conjugacy class of $RH(F(\cdot,1))$; cf. Corollary~\ref{cor:main2}.
\par
In a restricted form,
Theorem~\ref{thm:main} first appeared in \cite{msw}. A system generated
by a $T$-periodic non-autonomous ordinary differential equation was considered there and
it was assumed that $T>h$ and $T/h$ is rational.
It motivated a construction of an algorithm for rigorous numerical computing the index. 
The general version of the theorem given in the present paper
provides a theoretical background for extensions of the algorithm to other classes of 
ordinary differential equations.

\par 
This paper is organized as follows. In Section~\ref{sec:basic} we
provide definitions and results 
necessary to formulate and prove the assertions of the paper.  Theorem~\ref{thm:main} is stated in Section~\ref{sec:main} and
its proof is provided in Section~\ref{sec:proof}. The proof consists of several steps; its outline
is presented at the beginning of that section. For convenience of the reader, in Appendix~\ref{sec:app} we added proofs of
some results of \cite{msw} under the notation and settings of the present paper.

\section{Basic notation and terminology}
\label{sec:basic}
\subsection{Set theory, topology, and algebra}
The symbols $\mathbb R$, $\mathbb Z$, and $\mathbb N$ denote, respectively, the sets of real numbers, 
 the set of integers, and  the set on non-negative integers. 
For a map $f$, the inclusion
$f\subset X\times Y$ means that 
$f\colon D\to Y$ and $D$, the domain of $f$, is a subset of $X$
(the case $D=\emptyset$ is not excluded). If $A\subset X$, the restriction of $f$ to $A\times Y$
is denoted as $f|_A$.   
If $X=Y$, we set $f^0:=\id_X$, the identity on $X$, and if $n\in \mathbb N$, $n\geq 1$, by $f^n$ we denote the $n$-th iterate $f\circ\cdots\circ f$. 
\par
A pair $(X,A)$ of topological spaces is called a \emph{topological pair} if $A\subset X$ and the topology
of $A$ is induced from $X$. By a continuous map $f\colon (X,A)\to (Y,B)$ we mean a continuous $f\colon X\to Y$
such that $f(A)\subset B$.  
A topological pair 
$(X,A)$ is called \emph{compact} provided both $X$ and $A$ are compact spaces.
For a topological pair $(X,A)$ we define its \emph{quotient space} $X/A$ as
\[
X/A:=X\setminus A\cup \{\ast\},
\]
where $\ast$ is a point outside of $X$, and endow $X/A$ 
with 
 the strongest topology for which the \emph{quotient
 map} $q\colon X\to X/A$, given by 
 $q(x)=x$ for $x\in X\setminus A$ and $q(A)=\ast$,
 is continuous. In particular, $X/\emptyset=X\cup \{\ast\}$ has the direct sum topology. 
 We assume that $\ast$ is the same for all constructions of the quotient spaces related to the results
 of the present paper. It follows, in particular, that $\ast$ is a fixed point of the
 map $X/A\to Y/B$ generated by a continuous map $(X,A)\to (Y,B)$.
\par
By a \emph{graded vector space} we mean a direct sum $\bigoplus_{p\in\mathbb Z} V_p$ of vector spaces.
A graded linear map $\phi\colon V\to W$ of graded vector spaces $V$ and $W$
is defined as $\bigoplus_p\phi_p$, where $\phi_p\colon V_p\to W_p$. A graded
matrix $A=\bigoplus_pA_p$ is defined as the block-diagonal matrix generated by the matrices $A_p$. If the bases of $V$ and $W$ are fixed, we identify a matrix $A$ with the linear map
determined by $A$ in these bases. In the whole text we omit the symbol $\circ$ when composing linear
maps.

\subsection{The Leray reduction}
By an \emph{endomorphism} we mean a (graded) linear map $\alpha\colon V_\alpha\to V_\alpha$ 
of a (graded) vector space $V_\alpha$ into itself. If, moreover, it is an isomorphism, it is called
an \emph{automorphism}. Two automorphisms $\alpha$ and $\beta$ are \emph{conjugated}, denoted
$\alpha\cong \beta$, if 
there exists an isomorphism $\gamma$ such that $\beta\gamma=\gamma\alpha$.
The \emph{Leray reduction} $R$ is a transformation
of $\alpha$ onto an automorphism $R\alpha$ defined as follows. Set
\[
\operatorname{gker}(\alpha):=\bigcup_{n\in\mathbb N} \ker(\alpha^n),
\]
the \emph{generalized kernel} of $\alpha$ and
set $\overline V_\alpha:=V_\alpha/\operatorname{gker}(\alpha)$. Let
$\overline \alpha$ denote the induced endomorphism of $\overline V_\alpha$. Define
\[
RV_\alpha:=
\operatorname{gim}(\overline \alpha):=\bigcap_{n\in\mathbb N} \operatorname{im}(\overline \alpha^n),
\]
the \emph{generalized image} of $\overline \alpha$, and define $R\alpha$ as 
$\overline \alpha$ restricted to $RV_\alpha\to RV_\alpha$. 

\begin{prop}
\label{prop:leray}
Let $\alpha\colon V_\alpha\to V_\alpha$ and $\beta\colon V_\beta\to V_\beta$ be linear endomorphisms, let 
$\gamma\colon V_\alpha\to V_\beta$ and $\delta\colon V_\beta\to V_\alpha$  be linear maps, 
and let $r\in\mathbb N$ be such that the diagrams
\[
\begin{tikzcd}[column sep=large]
V_\alpha \ar[r,"\alpha"] \ar[d,"\gamma"]
&
V_\alpha \ar[d,"\gamma"]
\\
V_\beta \ar[r,"\beta"]
&
V_\beta
\end{tikzcd}
\qquad
\begin{tikzcd}[column sep=large]
V_\alpha \ar[r,"\alpha^r"] \ar[d,"\gamma"]
&
V_\alpha \ar[d,"\gamma"]
\\
V_\beta \ar[r,"\beta^r"] \ar[ru,"\delta"]
&
V_\beta
\end{tikzcd}
\]
commute. 
Then $R\alpha\cong R\beta$.
\end{prop}
Actually, the same result can be obtained for a general class of functors, called
\emph{normal consistent retractors} and a proof of that more general result is contained
in \cite[Appendix A]{msw}.  Below, in Subsection~\ref{ssec:a1} of the appendix we provide a direct
proof of Proposition~\ref{prop:leray}.
 
\subsection{Singular homologies}
\label{ssec:singular_h}
In the present paper $H$ denotes the singular homology functor with coefficients in a field $\mathbb F$.
Recall briefly its construction (we refer to \cite{dold} for a detailed exposition). 
Let $\Delta^d$ denote the standard $d$-dimensional simplex, i.e.
the convex envelope of the canonical basis of $\mathbb R^{d+1}$. A \emph{singular simplex} on a
topological space $X$ is a continuous map $\Delta^d\to X$. By 
$S(X)=\bigoplus_d S_d(X)$, the graded 
vector space of \emph{singular
chains}, we mean the set of all formal linear combinations of singular simplexes with coefficients in $\mathbb F$.
This is a chain complex and the corresponding boundary operator is denoted $\partial$.
The \emph{support} $|c|$ of
a chain $c=\sum_i\alpha_i\sigma_i\in S_d(X)$
is defined as
\[
|c|:=\bigcup_{i\colon\alpha_i\neq 0}\sigma_i(\Delta^d).
\]
A continuous map $f\colon X\to Y$ induces a chain map $S(f)\colon S(X)\to S(Y)$.
Let $A\subset X$ and let $i\colon A\hookrightarrow X$ be the inclusion map.
If $A\subset X$, $c\in S(X)$, and $|c|\subset A$, we treat $c$ also as 
an element of $S(A)$, 
hence $S(i)$ is equal to
the inclusion homomorphism $S(A)\hookrightarrow S(X)$. 
For the pair $(X,A)$, the vector spaces of
singular \emph{cycles}, \emph{boundaries}, and \emph{homologies} are denoted, respectively,
\begin{align*}
&
Z(X,A):=\{x\in S(X)\colon \partial z\in S(A)\},
\\
&
B(X,A):=\partial S(X)+S(A),
\\
&
H(X,A):=Z(X,A)/B(X,A).
\end{align*}
The \emph{homology class} of  $z\in Z(X,A)$ in $H(X,A)$ is denoted $[z]$.
A continuous map $f\colon (X,A)\to (Y,B)$ induces a linear map $H(f)\colon H(X,A)\to H(Y,B)$ 
given by $H(f)[z]:=[S(f)z]$. For $t\in [0,1]$, let $(\id,t)\colon (X,A)\to (X,A)\times [0,1]$ be defined
as $x\to (x,t)$.
 By a \emph{prism operator}
$P\colon S(X,A)\to S((X,A)\times [0,1])
$ 
we mean a chain homotopy joining the chain maps $S(\id,0)$ and $S(\id,1)$, i.e.
\[
S(\id,0)-S(\id,1)=\partial P+P\partial.
\]
 \par 
 A metrizable space $X$ is called an \emph{absolute neighborhood retract} (shortly: an ANR) if
 it is a neighborhood retract in every metrizable space in which it is embedded as a closed subset. 
 If $(X,A)$ is a pair of
 compact ANRs then its homologies are finite dimensional and 
 the quotient map induces an isomorphism $H(X,A)\to H(X/A,\ast)$.

\subsection{Local dynamical systems}
Let $X$ be a topological space. By  a \emph{local dynamical
system} (or, more precisely, a \emph{continuous-time local dynamical system}) on 
$X$ we mean a continuous map $\phi\colon D\to X$, where $D$ is an open
subset of $X\times \mathbb R$ and the set $D_x:=\{t\in\mathbb R\colon (x,t)\in D\}$ is
an open interval containing $0$ 
for each $x\in X$,
such that if $t\in D_x$ then $s+t\in D_x$ if and only
if $s\in D_{\phi(x,t)}$, and 
\begin{equation}
\label{eq:phi}
\phi(x,0)=x,\quad
 \phi(\phi(x,t),s)=\phi(x,s+t).
 \end{equation}
 If $D=X\times\mathbb R$, $\phi$ is called simply a \emph{dynamical system}.
By a \emph{semi-dynamical system} (or \emph{continuous-time semi-dynamical system}) on $X$ we mean a continuous map $\phi\colon X\times [0,\infty)\to X$ satisfying
\eqref{eq:phi}. 
For $A\times J\subset D$ we write $\phi(A,J)$ instead of $\phi(A\times J)$.
If $J$ is an interval contained in the domain of $\phi(x,\cdot)$,
the set $\phi(x,J)$ is called a \emph{segment}. 
Its \emph{length} is equal to the length of $J$. 
By $\phi(x)$ and $\phi^+(x)$, the 
\emph{trajectory} and, respectively, the \emph{positive semi-trajectory} of $x$, we mean the segment 
$\phi(x,J)$ with $J=D_x$ and, respectively, $J=D_x\cap [0,\infty)$,
i.e. the set of all $\phi(x,t)$ such that $\phi(x,\cdot)$ is defined at $t$ and, if the 
positive semi-trajectory is considered,
$t\geq 0$. In the case $D_x=\mathbb R$, $\phi(x)$ is called a \emph{full trajectory}.
\par 
For a continuous map $f\colon D\to X$, where $D$ is open in $X$,  the sequence $f^n$ of the iterates
of $f$ 
is called a 
\emph{discrete-time local semi-dynamical
system} (shortly: $f$ is a discrete-time local dynamical system).
If, moreover, the domain of $f$ is equal to $X$ (i.e. $f\colon X\to X$), it
 is called a \emph{discrete-time semi-dynamical system}. In this case
 a \emph{positive semi-trajectory} of $x$ is defined as the set of all $f^n(x)$, $n\in\mathbb N$,
 and
by a \emph{full trajectory} of $x$ we mean a set $\{x_n\colon n\in\mathbb Z\}$ such that 
$x_0=x$, and $f(x_n)=x_{n+1}$ for each $n\in\mathbb Z$.

\subsection{Isolated invariant sets, index pairs, and the Conley index}
\label{ssec:isolated}
In the sequel we assume that $X$ is a locally compact metrizable space $X$.
Let $S\subset X$. 
Both in the continuous-time and discrete-time dynamical systems, 
$S$ is called \emph{invariant} if each $x\in S$ has a full trajectory and that trajectory
is contained in $S$. 
An example of invariant set with respect to a continuous-time system $\phi\colon D\to X$
is the \emph{omega-limit set} of a point $x\in X$ for which
$[0,\infty)\subset D_x$. It is defined as  
\[
\omega(x):=\bigcap_{t>0}\clos\phi(x,[t,\infty)). 
\]
\par
Let $A\subset X$ and let $\phi$ be an arbitrary local dynamical system.
The \emph{invariant part} of $A$ is the maximal invariant set
contained in $A$. It is denoted as $\Inv_\phi A$ or $\Inv A$ if it does not lead to confusions.
If $A$ is closed then $\Inv A$ is also closed.
A compact set $N\subset X$ is called an \emph{isolating neighborhood}
if $\Inv N\subset \interior N$; equivalently:
no trajectory contained in $N$ intersects the boundary of $N$.
A compact set $S\subset X$ is called 
\emph{isolated invariant} if there exists an isolating neighborhood $N$
such that $S=\Inv N$. In this case $N$ is called an
\emph{isolating neighborhood of $S$}. 
\par
A pair $(N,L)$ of compact subsets of $X$ is called \emph{isolating} if
\[
\Inv\clos(N\setminus L)\subset \interior(N\setminus L).
\]
An isolating pair $(N,L)$ 
is called an \emph{index pair} if 
\emph{$L$ is positively invariant in $N$} and \emph{$L$ is an exit set of $N$}. In the 
case of a
continuous-time local dynamical system $\phi\colon D\to X$ 
this means
\[
\forall x\in L,\, t>0\colon \{\phi(x,[0,t])\subset N \Rightarrow \phi(x,t)\in L\}
\]
($L$ is positively invariant in $N$)
and
\[
\forall x\in N,\, t>0\colon  \{\phi(x,t)\notin N \Rightarrow \exists t'\in [0,t)\colon
\phi(x,[0,t'])\subset N,\  \phi(x,t')\in L\}
\]
($L$ is an exit set of $N$),
while in the case of a discrete-time local dynamical system $f\colon D\to X$ it is assumed that 
the set $N$ is contained in the domain $D$ and
\begin{align*}
&f(L)\cap N\subset L\quad
\text{($L$ is positively invariant in $N$),}&
\\
&\clos(f(N)\setminus N)\cap N\subset L 
\quad \text{($L$ is an exit set of $N$).}&
\end{align*}
(For a pair satisfying the above two conditions the term ``weak index pair'' was used in \cite{msw}.)
If $S=\Inv\clos(N\setminus L)$, the index pair $(N,L)$ is called an \emph{index pair for $S$}.
The basic feature of index pairs is recalled in the following two results.  

\begin{prop}
\label{prop:ip_phi}
If $(N,L)$ is an index pair with respect to a continuous-time local dynamical system $\phi$ then
the map 
\begin{align*}
&\phi_{(N,L)}\colon N/L\times [0,\infty)\to N/L,
\\
&
\phi_{(N,L)}(x,t):=\begin{cases}
\phi(x,t),&\text{if $\phi(x,[0,t])\subset N\setminus L$},
\\
\ast,&\text{otherwise}
\end{cases}
\end{align*}
is a semi-dynamical system on $N/L$.
\qed
\end{prop}
In the above statement 
``otherwise'' includes the case $(x,t)$ is not contained in the domain of $\phi$. 
Similarly through the rest of the paper. 
\begin{prop}
\label{prop:ip_f}
If $(N,L)$ is an index pair with respect to a discrete-time local dynamical system $f$ then 
the map
\begin{align*}
&f_{(N,L)}\colon N/L\to N/L,
\\
&
f_{(N,L)}(x):=\begin{cases}
f(x),&\text{if $x,f(x)\in N\setminus L$},
\\
\ast,&\text{otherwise}
\end{cases}
\end{align*}
is continuous.
\qed
\end{prop}

The \emph{homological Conley index} of an invariant set $S$ in a discrete-time local dynamical
system $f$, denoted $CH(S,f)$, is defined as the conjugacy class of the automorphism $RH(f_{(N,L)})$ for
an index pair $(N,L)$, where $R$ is the Leray reduction. 
It does not depend on the choice of $(N,L)$.

\subsection{Rotating systems and Poincar\'e maps}
\label{ssec:poincare}
We represent the circle as the quotient group $\mathbb R/\mathbb Z$.
Assume that $\theta\colon X\to \mathbb R/\mathbb Z$ is a continuous
map. For $a\in \mathbb R$ set $X_a:=\theta^{-1}(a+\mathbb Z)$,
hence $X_a=X_b$ if $a=b\bmod 1$. 
Define another locally compact metrizable space 
\[
\widetilde X:=\{(x,a)\in X\times \mathbb R\colon x\in X_a\}.
\]
Let 
\[
\zeta\colon \widetilde X\ni (x,a)\to x\in X,\quad
\widetilde \theta\colon \widetilde X\ni (x,a)\to a\in \mathbb R 
\]
be the projections.
The map $\zeta$ is a covering and the diagram
\[
\begin{tikzcd}
\widetilde X \ar[d,"\zeta"]
\ar[r,"\widetilde\theta"]
&
\mathbb R \ar[d,"(\cdot)+\mathbb Z"]
\\
X \ar[r,"\theta"]
&
\mathbb R/\mathbb Z
\end{tikzcd}
\]
 commutes.
Let $\phi\colon D\to X$ be a
continuous-time local dynamical system on $X$. 
It  
generates a local dynamical system $\widetilde\phi$ on $\widetilde X$
as follows.
Let $y\in \widetilde X$ and let $t\in D_{\zeta(y)}$. Assume $t>0$, hence $\{\zeta(y)\}\times [0,t]\subset D$. Set $\widetilde\phi(y,t)$ as
the terminal point of the lift of the path $[0,t]\ni s\to \phi(\zeta(y),s)\in X$
starting at $y$. 
In an analogous way we deal with the case $t<0$.
Thus $D_y$, the domain of $\widetilde\phi(y,\cdot)$, is equal to $D_{\zeta(y)}$.
In the sequel we assume that $\phi$ is \emph{rotating}, which means that
for every $y\in Y$ and $t>0$ such that $t\in D_{\zeta(y)}$,
\[
\widetilde\theta(y)<\widetilde\theta\widetilde\phi(y,t).
\]
\begin{prop}
\label{prop:omega_limit}
For every $y\in \widetilde X$, 
the omega-limit set $\omega(y)$ is empty. \qed
\end{prop}

For $a\in\mathbb R$ define $\widetilde X_a$ as $\widetilde\theta^{-1}(a)$, hence
$\widetilde X_a:=X_a\times a$. Let
$y\in \widetilde X$ and assume that $\widetilde\phi(y)\cap \widetilde X_a\neq \emptyset$. 
Define the \emph{access time} of $y$ to
$\widetilde X_a$ as
\[
\tau_a(y):=
\begin{cases}
\sup\{t\in D_y\colon \widetilde\theta\widetilde\phi(y,t)\leq a\}\in [0,\infty),&
\text{if $a\geq \w \theta(y)$,}
\\
\inf\{t\in D_y\colon \widetilde\theta\widetilde\phi(y,t)\geq a\}\in (-\infty,0],&
\text{if $a\leq \w\theta(y)$.}
\end{cases}
\]
i.e. $\tau_a(y)$ is equal to the time at which the trajectory of $y$ reaches $\widetilde X_a$.
Define a map
\[
\tau\subset (\widetilde X \times \mathbb R) \times \mathbb R,
\quad \tau(y,a):=\tau_a(y).
\]

\begin{prop}
\label{prop:tau_continuous}
The map $\tau$ is continuous, its domain is open, and the map
\[
\widetilde\Phi\subset (\widetilde X\times\mathbb R)\times \widetilde X,\quad \widetilde\Phi(y,t):=
\widetilde\phi(y,\tau(y,\w\theta(y)+t))
\]
is a local dynamical system on $\widetilde X$.
\qed
\end{prop}
Let $t$ be a real number.
 If $a=b\bmod 1$ then $\tau_{a+t}(x,a)=\tau_{b+t}(x,b)$, hence 
\[
T_t(x):=\tau_{a+t}(x,a)
\]
is independent of the choice of $a\in \theta(x)$.
By Proposition~\ref{prop:tau_continuous} the formula 
\[
\Phi(x,t):=\phi(x,T_t(x)).
\]
defines another rotating local dynamical system
$\Phi$ on $X$. We call it the \emph{translation system} associated to $\phi$ and we call 
$
\Phi_t:=\Phi(\cdot,t)  
$
a \emph{translation map}. 
Even if $\phi$ is a dynamical system, $\Phi$ can be local, i.e. its domain is an essential
subset of $X\times \mathbb R$.
For $a\in\mathbb R$ we treat the translation map as a map
\[
\Phi_t\subset X_a\times X_{a+t}.  
\]
In particular, if $t=1$ we call $T:=T_1$ the \emph{return time map}  and
\[
\Phi_1=\phi(\cdot,T(\cdot))\subset X_a\times X_a
\]  
a \emph{Poincar\'e map}. In the case $a=0$, we denote it by $\Pi$. Thus, for $n\in\mathbb N$,
\[
\Pi^n=\Phi_n\subset X_0\times X_0.
\] 
For $h\in \mathbb R$ we denote by $\phi^h$
the map $\phi(\cdot,h)\subset X\times X$. If $h>0$, we call it a \emph{discretization} of $\phi$.

\begin{prop}
\label{prop:s_invariant}
Let $S\subset X$. The following conditions are equivalent.
\begin{itemize}
\item[(a)] $S$ is an isolated invariant set with respect to $\phi$,
\item[(b)] for every (equivalently: for some) $h\neq 0$, $S$ is an isolated invariant set with respect to $\phi^h$,
\item[(c)] $S$ is an isolated invariant set with respect to 
$\Phi$,
\item[(d)] for every (equivalently: for some) $t\neq 0$, $S$ is an isolated invariant set with respect to $\Phi_t$,
\item[(e)] for every (equivalently: for some) $a\in\mathbb R$, $S_a$ is an isolated invariant set with respect to 
the Poincar\'e map $\Phi_1\subset X_a\times X_a$ and $S_{a+t}=\Phi_t(S_a)$
for each $t$.
\end{itemize}
\end{prop}
The equivalence of the conditions (a) and (b) was proved in \cite{mrozek-rm} (this result is valid for all local dynamical systems), hence also (c) holds if and only if (d) is satisfied.
Since the sets of trajectories of $\phi$ and $\Phi$ are equal each to the other, $S$ has the same
isolating neighborhoods in both systems and the equivalence of (a) and (c) is obvious. The equivalence
of (c), (d), and (e) is essentially stated as
 \cite[Prop.\,4.1]{msw}. 
We complete the proof of Proposition~\ref{prop:s_invariant} in
Subsection~\ref{ssec:s_inv} of the appendix.
\par
Since  the map $t\to T_t(x)$ is increasing for each $x\in X$,
the following result holds.
\begin{prop}
\label{prop:index_pairs_phi}
A pair 
$(N,L)$ of compact subsets of $X$ is an index pair with respect to $\phi$ if and only if it is an index pair with respect
to the translation system $\Phi$.\qed
\end{prop}

\subsection{Contiguous cycles}
\label{ssec:contiguous}
For $Z\subset X$ and $a\in\mathbb R$ define 
\[
Z_a:=Z\cap X_a,\quad \widetilde Z_a=Z_a\times a.
\]
Let $J\subset \mathbb R$ be an interval. Define
\[
\widetilde Z_J:=\{(x,a)\in \widetilde X\colon a\in J,\ x\in Z_a\},
\]
in particular, $\widetilde X_J=\widetilde\theta^{-1}(J)$.
In the case $J=\mathbb R$ we write $\widetilde Z$ instead of $\widetilde Z_{\mathbb R}$, hence
\[
\widetilde Z_J=\widetilde Z\cap \widetilde \theta^{-1}(J)=\widetilde Z\cap \widetilde X_J.
\]
Let $h>0$,
let $(N,L)$ be an index pair with respect to $\phi^h$, and let $a,b\in\mathbb R$, $a\leq b$. 
Assume that $u\in Z(N_a,L_a)$ and $v\in (N_b,L_b)$.
The pair $(u,v)$ is called a \emph{pair of 
contiguous cycles over $[a,b]$} if there exist chains $c\in S(\widetilde N_{[a,b]})$
and $d\in S(\widetilde L_{[a,b]})$ such that
\begin{equation}
\label{eq:contiguous_cycles}
u\times a-v\times b=\partial c+d,
\end{equation}
where $u\times a\in Z(\widetilde N_a,\widetilde L_a)$ 
and $v\times b\in Z(\widetilde N_b,\widetilde L_b)$ are induced by the maps $x\to (x,a)$ and,
respectively,
$x\to (x,b)$.
A pair of contiguous cycles $(u,v)$ is \emph{$h$-movable} if for some chains $c$ and $d$
satisfying \eqref{eq:contiguous_cycles},
\[
\widetilde\phi(|c|,[0,h])\subset \widetilde N,\quad
\widetilde \phi(|d|,[0,h])\subset \widetilde L.
\]
We call such chains $c$ and $d$ \emph{associated} to $(u,v)$.

\begin{prop}
\label{prop:movable}
If $(u,v)$ is an $h$-movable pair of contiguous cycles over $[a,b]$ 
and $a<b$ then 
\begin{align*}
&\widetilde\phi(|u\times a|,[0,h])\subset \widetilde N,\quad
\widetilde\phi(|\partial u\times a|,[0,h])\subset \widetilde L,
\\
&
\widetilde\phi(|v\times b|,[0,h])\subset \widetilde N,\quad
\widetilde\phi(|\partial v\times b|,[0,h])\subset \widetilde L.
\end{align*}
\end{prop}
\begin{proof}
Since
\[
\widetilde\phi(|u\times a-v\times b|,[0,h])=\widetilde\phi (|\partial c-d|,[0,h])\subset 
\widetilde\phi(|\partial c|,[0,h])\cup \widetilde\phi(|d|,[0,h])\subset
\widetilde N
\]
and $|u\times a|\cap |v\times b|=\emptyset$, two of the required inclusions are satisfied. Since
\[
\widetilde\phi(|\partial u\times a-\partial v\times b|,[0,h])=\widetilde\phi (|\partial d|,[0,h])\subset \widetilde L
\]
and $|\partial u\times a|\cap |\partial v\times b|=\emptyset$, the other two inclusions also hold.
\end{proof}

The following two results are also consequences of the fact that
the support of the sum of chains is contained in the union of their supports.

\begin{prop}
\label{prop:contiguous_abc} 
If $(u_{i-1},u_i)$ is 
an $h$-movable pair of contiguous cycles over $[a_{i-1},a_i]$, where $a_{i-1}\leq a_i$
for $i=1,\ldots ,n$,
then $(u_0,u_n)$ is an $h$-movable pair of contiguous cycles over $[a_0,a_n]$.\qed
\end{prop}

\begin{prop}
\label{prop:sum_contiguous}
If $(u_i,v_i)$ is an $h$-movable pair of contiguous cycles and $\lambda_i\in\mathbb F$ for
$i=1,\ldots,n$ then $(\sum_i\lambda_iu_i,\sum_i\lambda_iv_i)$ is an $h$-movable pair of contiguous
cycles.\qed
\end{prop}

\section{The main theorem}
\label{sec:main}
Recall that $\phi$ is a
rotating local dynamical system on $X$.
It follows by Proposition~\ref{prop:s_invariant},  that if $S$ is an isolated invariant set for $\phi^h$ with some $h>0$ then 
$S_0$ is an isolated invariant set for the Poincar\'e map $\Pi$.

\begin{theorem}
\label{thm:main}
Let
$h>0$ and let $(N,L)$ be an index pair for an isolated invariant
set $S$ with respect to $\phi^h$. If $N_0$ and $L_0$ are ANRs, 
$
n=\dim H(N_0,L_0),
$
$A=[a_{ij}]$ is a graded $(n\times n)$-matrix over $\mathbb F$, and
$
\left(u_j,\sum_{i=1}^n a_{ij}u_i\right)
$
for $j=1,\ldots,n$
is an $h$-movable pair of contiguous cycles over $[0,1]$ such that
$\{[u_j]\colon j=1,\dots,n\}$ is a basis of $H(N_0,L_0)$ 
then $CH(S_0,\Pi)$ is equal to the conjugacy class of the Leray reduction $RA$.
\end{theorem} 
Under the assumptions that
$X=X_0\times \mathbb R/\mathbb Z$, the return time map $T$ is 
constant, $T>h$, and $T/h$ is a rational number, a proof of Theorem~\ref{thm:main} is essentially given in \cite{msw}. In Section~\ref{sec:proof}
we prove the theorem in full generality. 

\begin{corollary}
\label{cor:main2}
Let $h>0$, $(N,L)$, and $(N_0,L_0)$ satisfy the assumptions
of Theorem~\ref{thm:main} and let $F_a\colon (N_0,L_0)\to (N_a,L_a)$
for $a\in [0,1]$ be a family of maps
such that
\[
F\colon (N_0,L_0)\times [0,1]\to (N,L),\quad F(x,a):=F_a(x)
\]
is a continuous. If $F_0=\id_{N_0}$ and 
\[
\phi^t(F_a(N_0))\subset N,\quad \phi^t(F_a(L_0))\subset L
\] 
for all $a\in [0,1]$ and $t\in [0,h]$ then
$CH(S_0,\Pi)$ is equal to the conjugacy class of $RH(F_1)$.
\end{corollary}
\begin{proof}
Let
\[
\widetilde F\colon (\widetilde N_0,\widetilde L_0)\times [0,1]\to (\widetilde N,\widetilde L)
\]
be the lift of $F$ and
let $u_i\in Z(N_0,L_0)$, $i=1,\ldots,n$ be such that $\{[u_i]\colon i=1,\ldots,n\}$ form a basis of $H(N_0,L_0)$. 
Then $(u_i,S(F_1)u_i)$ is an $h$-movable pair of contiguous cycles over $[0,1]$, since 
we can set $S(\widetilde F)P(u_i\times 0)\in S(\widetilde N_{[0,1]})$ and 
$S(\widetilde F)P\partial(u_i\times 0)\in S(\widetilde L_{[0,1]})$ as associated chains, where
\[
P\colon S(\widetilde N_0,\widetilde L_0)\to S((\widetilde N_0,\widetilde L_0)\times [0,1])
\]
is a prism operator.
Let $A$ be the matrix of $H(F_1)$ with respect to the basis $\{[u_i]\}$. Then there exist chains
$c\in S(\widetilde N_1)$ and $b\in S(\widetilde L_1)$ such that
\[
S(F_1)u_i\times 1-Au_i\times 1=\partial c+b,
\]
hence $(u_i,Au_i)$ is also an $h$-movable pair of contiguous cycles over $[0,1]$ 
by
Proposition~\ref{prop:contiguous_abc} and
the assumed inclusions $\phi^t(N_0)\subset N$ and $\phi^t(L_0)\subset L$ for $t\in [0,h]$. 
The result thus follows from Theorem~\ref{thm:main}.
\end{proof}

\section{Proof of Theorem~\ref{thm:main}}
\label{sec:proof}

\subsection{An outline of the proof}
Here we sketch the main steps of the proof; sometimes they appear in a different 
form and in a different order in the rigorous arguments presented in the 
following subsections.
\par\noindent
\textbf{1st step.} We should find an index pair for $\Pi$. For this purpose
we construct sets $N'$, $L'$, $N''$, and $L''$ such that 
\[
(N,L)\supset (N',L) \subset (N',L') \supset (N'',L'')
\]
and
$(N',L)$ is an index pairs with respect to $\phi^h$,
$(N',L')$ is an index pair with respect to both $\phi$ and $\phi^h$,
and $(N'',L'')$ is an index pair for the translation map $\Phi_1$. 
Then $(N''_0,L''_0)$ is the required index pair. By the construction, $N'_0/L'_0=N''_0/L''_0$, 
hence  $CH(S_0,\Pi)$ is equal to the conjugacy class
of $RH(\Phi'_1)$, where $\Phi'_1$ is the self-map of $N'_0/L'_0$ induced by 
$\Phi_1$.   
\par\noindent
\textbf{2nd step.} In order 
prove that $RH(\Phi'_1)\cong RA$ we
apply Proposition~\ref{prop:leray}. We should 
determine an
$r\in\mathbb N$ and homomorphisms $\gamma$ and $\delta$ between
$H(N_0,L_0)$ and $H(N'_0/L'_0,\ast)$ which appear in the corresponding
commutative diagrams.
We set $r=k+m$, where $k$ and $m$ relate to the maximal sizes of the sequences 
$x,\phi^h(x),\phi^{2h}(x),\ldots$
 contained in $L'\setminus L$ and $N\setminus N'$, respectively.
 The map $\gamma$ is  constructed as follows. Denote by $U$ the 
 linear span of the cycles $u_i$; we treat the matrix $A$ as an endomorphisms of $U$
 generated by the map $u_i\to \sum_j a_{ij} u_i$.
For $u\in U$ and $\ell\in \mathbb N$, the pair of contiguous cycles $(u,A^\ell u)$
over $[0,\ell]$ is $h$-movable. It follows, in particular, that
$u\in U$ can be treated as a chain in $N'$ and therefore also as a cycle in $(N'_0,L_0)$. 
Define the class $\gamma[u_i]$ as the class of 
 the cycle induced from $u_i$ by the quotient map $N'\to N'/L'$. The classes
 $[u_i]$ form a basis of $H(N_0,L_0)$, 
hence the homomorphism $\gamma$ is defined.
The map $\delta$ will be determined in the 6th step.
\par\noindent
\textbf{3rd step.}
Since we will use $h$-movable contiguous pairs of cycles, one should 
lift the considered index pairs in the space $X$ to the corresponding pairs in $\widetilde X$.
These pairs are non-compact, hence in order to keep the notation
of the theory of    
isolated invariant sets it is convenient to switch to the quotient space $Y$ of $\widetilde X$ 
corresponding to the quotient homomorphism $\mathbb R\to \mathbb R/p\mathbb Z$ for
$p>k+m$ large enough. There are derived systems $\psi$ and $\Psi$ in $Y$ corresponding 
to $\phi$ and $\Phi$, and the index pairs $(Q,P)$, $(Q',P)$, etc. corresponding to
the index pairs $(N,L)$, $(N',L)$, etc., respectively. All dynamics related to $\psi$ within the range
of $Y$ over the interval $[0,p)$ is the same as the corresponding dynamics of $\widetilde \phi$.
In particular, the notion of the $h$-movability of contiguous cycles over intervals 
contained in $[0,p)$ has an equivalent formulation in terms of 
$Y$ and $\psi$.
\par\noindent
\textbf{4th step.} For $\ell=0,\ldots,k+m$ we define homomorphisms $\Gamma_\ell$
from $H(Q_\ell,P_\ell)$ to $H(Q'_\ell/P'_\ell,\ast)$ as the composition of   
$\gamma\colon H(N_0,L_0)\to H(N_0/L_0,\ast)$ with the isomorphisms  
induced by the map $x\to (x,\ell)$ and its inverse,
and for $\ell<k+m$ we denote also by $A$ a homomorphisms $H(Q_\ell,P_\ell)\to H(Q_{\ell+1},P_{\ell+1})$ 
corresponding to the endomorphism of $H(N_0,L_0)$ induced by the matrix $A$. In order to get
the equation $H(\Phi'_1)\gamma=\gamma A$
we should
prove that $H(\Psi'_1)\Gamma_0=\Gamma_1A$.
To this purpose we denote by $C$ and $B$ the sets
of points $y$ such that $\psi(y,[0,h])\subset Q$ and $\psi(y,[0,h])\subset P$, respectively, and we define
a map $F$ from $(C_{[0,1]},B_{[0,1]})$ to 
$(Q'_1/P'_1,\ast)$ which collapses $C_{[0,1]}$ along the trajectories
of the semi-dynamical system $\psi_{(Q',P')}$. Then the homomorphism 
$S(F)$ applied to chains associated 
to each pair $(u_i,Au_i)$ of contiguous cycles  over $[0,1]$ provides the required equality.
\par\noindent
\textbf{5th step.}
In order to define the map $\delta$ required in the 2nd step, at first
we construct a continuous map $G$ from $(C_{[0,k]},B_{[0,k]})$ to $(Q'_k/P_k,\ast)$.
If  $\tau_k$ is defined at $y\in C_{[0,k]}$, there are $\overline s_y\in [0,h)$ and
$\overline \ell\in\mathbb N$ such that $\tau_k(y)=\overline s_y+\overline \ell_yh$. 
We set $G(y)$
as the image of $\psi(y,\overline s_y)$ under the $\overline\ell_y$-th iterate 
the map $\psi^h_{(Q',P)}$. The number $k$ and the set $P'$ are chosen in the way that
$G|_{C_0}$ factors through $C_0/P'_0=Q'_0/P'_0$. The induced continuous 
map $Q'_0/P'_0\to Q'_k/P_k$ is denoted by $r$.
We prove directly that $r$ composed with the quotient map $j'\colon Q'_k/P_k\to Q'/P'_k$ is equal to $\Psi'_k$.
For $u\in U$, we denote by $\overline u^k$ the chain induced from $u\times k$ by
the quotient map $Q'_k\to Q'_k/P_k$. We prove that the classes $[\overline u_i^k]$ are linearly independent
in $H(Q'_k/P_k,\ast)$.
In an analogous way as $S(F)$ in the 4th step, the homomorphism
$S(G)$ applied to chains associated to the pair of contiguous cycles $(u_i,A^ku_i)$ over $[0,k]$ 
provides the equation $H(r)\Gamma_0[u_i\times 0]=[\overline v_i^k]$ in $H(Q'_k/P_k,\ast)$, where
$v_i=A^ku_i$. 
\par\noindent
\textbf{6th step.} We take an arbitrary homology class $[z]$ in $H(Q'_k/P_k,\ast)$. Its image
under the homomorphism induced by the inclusion map $Q'_k/P_k\to Q_k/P_k$ is equal
to the image of $[\overline u^k]$ for exactly one $u\in U$. 
We define the map
$\Delta\colon H(Q'_0/P'_0,\ast)\to H(Q_{k+m},P_{k+m})$ as the composition
of $H(r)$, the map $[z]\to [u\times k]$, and $A^m$.
Consequently, we get the equation $\Delta\Gamma_0=A^{k+m}$ by the 5th step.
By the choice of $m$ we are able
to construct a continuous map $d\colon Q_k/P_k\to Q'_{[k,k+m]}/P_{[k,k+m]}$ such
that $S(d)z$ is homologous to $S(d)\overline u^k$. Furthermore, there is a continuous
map $J$ which collapses $Q'_{[k,k+m]}/P'_{[k,k+m]}$ to $Q'_{k+m}/P'_{k+m}$ along
the trajectories of the semi-dynamical system $\psi_{(Q',P')}$.
We 
apply $S(J)$ to the above homologous cycles and as a result we conclude that 
$H(\Psi'_m)\Gamma_k[u\times k]=H(\Psi'_m)H(j')[z]$.
Since $j'\circ r=\Psi'_k$ and $H(\Psi'_m)\Gamma_k=\Gamma_{k+m}A^m$
by the 4th step, the equation $\Gamma_{k+m}\Delta=H(\Psi'_{k+m})$
holds. 
Thus $\delta$, defined as $\Delta$ composed with
 the isomorphisms induced by the maps $x\to (x,0)$ and $(x,k+m)\to x$, 
satisfies $\delta\gamma=A^{k+m}$ and $\gamma\delta=H(\Phi'_1)^{k+m}$, 
hence the theorem follows.

\subsection{Auxiliary index pairs}
Define 
\begin{align*}
&K:=\bigcup\{\sigma\subset N\colon \text{$\sigma$ is a segment starting at a point in $L$}\},
\\
&M:=\bigcup\{\sigma\subset N\colon\text{$\sigma$ is a segment of length $h$}\}.
\\
&N':=K\cup M.
\end{align*} 

\begin{prop}
\label{prop:index_pairs}\
\begin{itemize}
\item[(a)] $(N',L)$ is an index pair for $S$ with respect to $\phi^h$,
\item[(b)] $\phi^h(N')\cap N\subset N'$,
\item[(c)] $(N',K)$ is an index pair for $S$ with respect to $\phi$.
\end{itemize}
\end{prop}
This result was proved in \cite[Sec.\,5.2]{msw}.
We provide its proof within the settings of the present paper in Subsection~\ref{ssec:a3} of the appendix.

\begin{lemma}\emph{(\cite{salamon}, Lemma~5.3)}
\label{lem:lyapunov}
There exists a continuous function $\alpha\colon N'\to [0,1]$ such that 
\begin{itemize}
\item[]
$\alpha(x)=1$ if and only if $\phi(x,[0,\infty))\subset N'$ and $\omega(x)\subset S$,
\item[]
$\alpha(x)=0$ if and only if $x\in K$,
\item[]
if $t>0$, $0<\alpha(x)<1$, and $\phi(x,[0,t])\subset N'$ then 
$\alpha(\phi^t(x))<\alpha(x)$.
\end{itemize}
\end{lemma}
Let  $\alpha$ be a function satisfying Lemma~\ref{lem:lyapunov}.
It follows, in particular, that 
the return time map $T$ is defined on the set $\alpha^{-1}(1)$. 
Let $W$ be a compact neighborhood
of $\alpha^{-1}(1)$ contained in the domain of $T$. Set 
\[
\lambda:=\max\{h,\max\{T(x)\colon x\in W\}\}.
\]
Choose a number $c<1$ such that 
if $\alpha(x)\geq c$ then $x\in W$ and
$
\phi(x,[0,\lambda])\subset N'\setminus K.
$
Set
\[
L':=\alpha^{-1}([0,c]),\quad
N'':=\alpha^{-1}([c,1]),\quad L'':=\alpha^{-1}(c),
\]
hence $T$ is defined on $\clos(N'\setminus L')$ and
\begin{equation}
\label{eq:phi_n-l}
\phi(\clos(N'\setminus L'),[0,\lambda])\subset N'\setminus K. 
\end{equation}

\begin{prop}
\label{prop:a4}
$(N',L')$ is an index pair for $S$ with respect to both $\phi$ and $\phi^h$.
\end{prop}
Again, essentially the same result was given in \cite[Sec.\,5.2]{msw}. We provide 
a proof of Proposition~\ref{prop:a4} in Subsection~\ref{ssec:a4}.
\par
As a direct consequence of the definition of $N''$ and $L''$, and Proposition~\ref{prop:index_pairs_phi},
one can verify that $(N'',L'')$ is an index pair for $S$ with respect to
 $\Phi_1$ (actually, also with respect to $\phi$, $\phi^h$, and $\Phi$). 
In particular, the requirement that $N''$ is contained in the domain of $\Phi_1$ is 
satisfied by \eqref{eq:phi_n-l}.  Consequently, an index pair for
the Poincar\'e map $\Pi$ is given in the following result.

\begin{prop}
\label{prop:ip_pi} 
$(N''_0,L''_0)$ is an index pair for $S_0$ with respect to $\Pi$.\qed
\end{prop}

Propositions~\ref{prop:index_pairs_phi} and \ref{prop:a4} imply that $(N',L')$ is 
an index pair for $\Phi$. Set
\[
\Phi':=\Phi_{(N',L')},
\]
hence 
for $a,t\in\mathbb R$ the map
$\Phi'_t=\Phi'(\cdot,t)\colon N'_a/L'_a\to N'_{a+t}/L'_{a+t}$ is given as
\[
\Phi'_t(x)=
\begin{cases}
\phi_{(N',L')}(x,T_t(x)),& \text{if $x\in N'_a\setminus L'_a$ and $T_t$ is defined at $x$},
\\
\ast,& \text{otherwise}
\end{cases}
\]
Since $N''/L''=N'/L'$,  
\[
\Pi_{(N''_0,L''_0)}=\Phi'_1\colon N'_0/L'_0\to N'_0/L'_0.
\]
Therefore, by Proposition~\ref{prop:ip_pi}, in order to prove Theorem~\ref{thm:main} it suffices to show
that 
\begin{equation}
\label{eq:required}
RH(\Phi'_1)\cong RA.
\end{equation}

\subsection{Reformulations of Theorem~\ref{thm:main}}
Set $f:=\phi^h$.

\begin{lemma}
\label{lem:k+m_nl}
There exist $\kappa$ and $\mu$ in $\mathbb N$ such that
\begin{align*}
&
\{x\in L'\setminus L\colon f(x),\dots,f^\kappa (x)\in N'\setminus L\}=\emptyset,
\\
&
\{x\in N\setminus N'\colon f(x),\ldots,f^\mu (x)\in N\setminus N'\}=\emptyset. 
\end{align*}
\end{lemma}
\begin{proof}
If for every $n$ there is a $x_n\in L'\setminus L$ such that 
$f(x_n),\ldots,f^n(x_n)\in N'\setminus L$, hence, by the positive invariance of $L'$ in $N$,
also $f(x_n),\ldots,f^n(x_n)\in L'\setminus L$, then $\omega(x)\subset \clos(L'\setminus L)$ for
each accumulation point $x$ of $\{x_n\}$, contrary to the fact that both $(N',L)$ and $(N',L')$ are isolating pairs 
for $S$, hence the first equation holds. 
Since $S\cap \clos(N\setminus N')=\emptyset$, 
the second statement follows by a similar argument.
\end{proof}
We claim that the return time map $T$ satisfies
\begin{equation}
\label{eq:inf_tau}
\inf\{T(x)\colon x\in N\}>0.
\end{equation}
Indeed,
assume $x_n\in N$ and $T(x_n)\to \inf\{T(x)\colon x\in N\}$.  Let $x_0\in N$ be
the limit of $x_n$.  If $T$ is defined at $x_0$ then there is nothing to prove. 
Assume the opposite.  Since $\phi(\cdot,h)$ is defined at $x_0$
(as it is required for an index pair with respect to $f$), there is a neighborhood $W$ of $(x_0,0)$ 
in $\widetilde X$ such that 
\[
\widetilde \phi(W,[0,h])\cap \widetilde X_1=\emptyset,
\]
hence
$T(x_n)> h$ for $n$ large enough and the inequality follows.
\par
We fix $\kappa$ and $\mu$ satisfying Lemma~\ref{lem:k+m_nl}.
Thanks to \eqref{eq:inf_tau}, there is a
$\nu\in\mathbb N$ such that 
\[
h<\nu \inf\{T(x)\colon x\in N\}.
\]
Define
\[
k:=\nu\kappa,\quad m:=\nu\mu.
\]
Since $H(\Phi'_1)^\ell=H(\Phi'_\ell)$ for every positive integer $\ell$, 
Proposition~\ref{prop:leray} implies that
the equation \eqref{eq:required}
 is proved if we show that  
the following result holds true.

\begin{prop}
\label{prop:main_nl}
There are homomorphisms
\[
\gamma\colon H(N_0,L_0)\to H(N'_0/L'_0,\ast),\quad
\delta\colon H(N'_0/L'_0,\ast)\to H(N_0,L_0)
\]
such that the diagrams
\[
\begin{tikzcd}[column sep=large,row sep=large]
H(N_0,L_0) \ar [r,"A"] \ar[d,"\gamma"']
&
H(N_0,L_0) \ar[d,"\gamma"]
\\
H(N'_0/L'_0,\ast) \ar[r,"H(\Phi'_1)"] 
&
H(N'_0/L'_0,\ast)
\end{tikzcd}
\
\begin{tikzcd}[column sep=large,row sep=large]
H(N_0,L_0) \ar [r,"A^{k+m}"] \ar[d,"\gamma"'] 
&
H(N_0,L_0) \ar[d,"\gamma"]
\\
H(N'_0/L'_0,\ast) \ar[r,"H(\Phi'_{k+m})"]
\ar[ru,"\delta"] 
&
H(N'_0/L'_0,\ast)
\end{tikzcd}
\]
commute.
\end{prop}
Let us fix  an integer $p\geq k+m+\nu$. We identify the space $\mathbb R/p\mathbb Z$
with the interval $[0,p)$ endowed with the 
corresponding topology and we construct a quotient
space $Y$ of 
$\widetilde X$ by identifying $(x,a)$ with $(x,b)$ if $a=b\bmod p$, i.e.
\[
Y=\{(x,a)\in X\times [0,p)\colon x\in X_a\}\subset X\times \mathbb R/p\mathbb Z.
\]
The maps $a\to a+p\mathbb Z\to a+\mathbb Z$ generate coverings
$\widetilde \zeta_p\colon \widetilde X\to Y$ and $\zeta_p\colon Y\to X$. 
Let $\psi$ be the local dynamical dynamical system on $Y$ defined
by lifting $\phi$ with respect to $\zeta_p$ in an analogous way as 
$\widetilde \phi$ was constructed.  This is a rotating system. Let $\Psi$ denote the translation
system associated to $\psi$; actually, $\Psi$ is induced on $Y$ from the system
$\widetilde \Phi$ given in Proposition~\ref{prop:tau_continuous}.  
\par
Let $a\in [0,p)$. 
Similarly as in the case of $\widetilde X$, we set $Y_a:=X_a\times a$ 
and, more generally, for $Z\subset Y$ and a closed interval $I\subset [0,a)$, 
\[
Z_I:=\{(x,a)\in Y\colon a\in I\}.
\]
The map $\widetilde \zeta_p$ 
identifies the set $\widetilde X_I$ with $Y_I$ as topological spaces and therefore
we do not distinguish between
the dynamics of $\widetilde \phi$
and $\psi$ (as well as $\widetilde \Phi$ and $\Psi$) restricted to these sets. 
In particular, Proposition~\ref{prop:omega_limit} has
the following reformulation.
\begin{lemma}
\label{lem:no_omega}
If $I$ is a closed interval in $[0,p)$ then
the system $\psi$ has no omega-limit sets contained in $Y_I$.\qed
\end{lemma}
We 
adapt from $\widetilde \phi$ the notion of an $h$-movable pair of contiguous cycles 
over any interval contained in $[0,k+m]$ and therefore all results stated below for $\psi$ are
valid also for $\widetilde \phi$. 
Define $g:=\psi^h$ and
\begin{align*}
&
S^\ast:=\{(x,a)\in Y\colon x\in S_a\},
\\
&Q:=\{(x,a)\in Y\colon x\in N_a\},\quad P:=\{(x,a)\in Y\colon x\in L_a\},
\\
&Q':=\{(x,a)\in Y\colon x\in N'_a\},\quad P':=\{(x,a)\in Y\colon x\in L'_a\},
\end{align*}
The set $S^\ast$ is isolated invariant in $Y$ and $(Q,P)$, $(Q',P)$,  and $(Q',P')$
are index pairs for $S^\ast$ with respect to $g$ by Propositions~\ref{prop:index_pairs}\,(a)
and \ref{prop:a4},
and $(Q',P')$ is an index pair for $S^\ast$ with respect to $\psi$ and $\Psi$
by Propositions~\ref{prop:index_pairs_phi} and \ref{prop:a4}.
Moreover, by Proposition~\ref{prop:index_pairs}\,(b),
\begin{equation}
\label{eq:q'q}
g(Q')\cap Q\subset Q'
\end{equation}
and
by \eqref{eq:phi_n-l} and the choice of $\kappa$ and $\mu$, 
\begin{align}
\label{eq:psi_q-p}
&
\psi(\clos(Q'\setminus P'),[0,h])\subset Q'\setminus P,
\\
\label{eq:choice_k}
&
\{y\in P'\setminus P\colon g(y),\dots,g^\kappa(y)\in Q'\setminus P\}=\emptyset,
\\
&
\label{eq:choice_m}
\{y\in Q\setminus Q'\colon g(y),\ldots,g^\mu(y)\in Q\setminus Q'\}=\emptyset.
\end{align}
The system $\psi$ and the map $g$ generate the semi-dynamical systems $g_{(Q,P)}$, $g_{(Q',P)}$,
$g_{(Q',P')}$, $\psi_{(Q',P')}$, and $\Psi':=\Psi_{(Q',P')}$. 
For $a,a+t\in [0,p)$ set
\[
\Psi'_t:=\Psi'(\cdot,t)\colon Q'_a/P'_a\to Q'_{a+t}/P'_{a+t},
\]
 hence
\[
\Psi'_t(y)=
\begin{cases}
\psi_{(Q',P')}(y,\tau_{a+t}(y)),& 
\text{if $y\in Q'_a\setminus P'_a$ and $\tau_{a+t}$ is defined at $y$},
\\
\ast,&\text{otherwise}.
\end{cases}
\]
For $\ell=0,\ldots,k+m-1$ let
$A$ denote the linear map $H(Q_\ell,P_\ell)\to H(Q_{\ell+1},P_{\ell+1})$ having the matrix 
$A$ in the 
bases $\{[u_i\times \ell]\}$ of $H(Q_\ell,P_\ell)$ and $\{[u_i\times (\ell+1)]\}$ of $H(Q_{\ell+1},P_{\ell+1})$.  
Using homomorphisms in homologies induced by the covering map $\zeta_p\colon Y\to X$ we reformulate
Proposition~\ref{prop:main_nl} to the following equivalent statement.
\begin{prop}
\label{prop:main_qp} 
For $\ell=0,\ldots, k+m$ there is a homomorphism 
\[
\Gamma_\ell\colon H(Q_\ell,P_\ell)\to H(Q_\ell,P_\ell)
\]
such that
\begin{itemize}
\item[(a)]
the diagram 
\[
\begin{tikzcd}[column sep=large,row sep=large]
H(Q_0,P_0) \ar [r,"A"] \ar[d,"\Gamma_0"']
&
H(Q_1,P_1) \ar[d,"\Gamma_1"]
\\
H(Q'_0/P'_0,\ast) \ar[r,"H(\Psi'_1)"] 
&
H(Q'_1/P'_1,\ast)
\end{tikzcd}
\]
commutes, 
\item[(b)]
there is a homomorphism $\Delta$
such that the diagram
\[
\begin{tikzcd}[column sep=large,row sep=large]
H(Q_0,P_0) \ar [r,"A^{k+m}"] \ar[d,"\Gamma_0"'] 
&
H(Q_{k+m},P_{k+m}) \ar[d,"\Gamma_{k+m}"]
\\
H(Q'_0/P'_0,\ast) \ar[r,"H(\Psi'_{k+m})"]
\ar[ru,"\Delta"] 
&
H(Q'_{k+m}/P'_{k+m},\ast)
\end{tikzcd}
\]
commutes.
\end{itemize}
\end{prop}
A proof of Proposition~\ref{prop:main_qp} (hence, as a consequence, also of Theorem~\ref{thm:main}) is
presented in the subsequent subsections.

\subsection{Beginning of the proof of Proposition~\ref{prop:main_qp}}

The identities and the inclusion map $e\colon Q'\hookrightarrow Q$ 
induce the following continuous maps in the quotient spaces: 
\begin{align*}
&j\colon Q'/P\to Q/P,\quad q\colon Q\to Q/P,
\\
&j'\colon Q'/P\to Q'/P',\quad q'\colon Q'\to Q'/P.
\end{align*}
It follows that $q'$ is a restriction of $q$ and $j'\circ q'$ is equal
to the quotient map $Q'\to Q'/P'$. 
(Actually, $j$ is an inclusion map since $Q'/P\subset Q/P$ and the quotient topology on $Q'/P$
 coincides with the topology induced from $Q/P$, but we will not take advantage
 of this fact later.)
In order to simplify notation, by the same 
letters $e$, $j$, $j'$, $q$, and $q'$ we denote also their restrictions, for example for $a\in [0,p)$ we
have
$e\colon Q'_a\hookrightarrow Q_a$,
$j\colon Q'_a/P_a\hookrightarrow Q_a/P_a$, 
$q\colon Q_a\to Q_a/P_a$, etc. 
In the sequel we will use also the sets
\[
C:=\{y\in Q\colon \psi(y,[0,h])\subset Q\},\quad B:=\{y\in Q\colon \psi(y,[0,h])\subset P\}.
\]
Clearly, $(C,B)\subset (Q',P)$.
\par
Set $U:=\left\{\sum_i\alpha_iu_i\colon \alpha_i\in\mathbb F\right\}$. This is the linear subspace
of $S(N_0)$ spanned on the cycles $u_i$,
$i=1,\ldots,n$, and $A$ is an endomorphism $U\to U$.
By assumption, the pair $(u_i,Au_i)$ is $h$-movable for each $u_i$,
hence 
for every $u\in U$ and all $\ell<\ell'\in \mathbb N$,
$(u,A^{\ell'-\ell} u)$ is an $h$-movable pair of contiguous cycles over $[\ell,\ell']$
as a consequence of Propositions~\ref{prop:contiguous_abc} and \ref{prop:sum_contiguous}.
In particular, by Proposition~\ref{prop:movable} 
for each $u\in U$ and $\ell=0,\ldots,k+m$,
\[
u\times\ell\in Z(C_\ell,B_\ell).
\]
A comment on the notation we use: in the sequel the chain $u\times \ell$ will appear as a cycle in 
$Z(C_\ell,B_\ell)$ or 
$Z(Q_\ell,P_\ell)$, or 
$Z(Q'_\ell,P_\ell)$. It will be clear from the context
in which of the homologies $H(C_\ell,B_\ell)$, $H(Q_\ell,P_\ell)$, and $H(Q'_\ell,P_\ell)$ 
the class 
$[u\times \ell]$ is treated.
\par
Since the map $x\to (x,\ell)$ is a homeomorphism $(N_0,L_0)\to (Q_\ell,P_\ell)$,
the classes $[u_i\times\ell]$ form a basis of $H(Q_\ell,P_\ell)$.
Moreover, since $(N_0,L_0)$
is a pair of compact ANRs, for $\ell=0,\ldots,k+m$ the map 
\[
H(q)\colon H(Q_\ell,P_\ell)\to H(Q_\ell/P_\ell,\ast)
\]
is an isomorphism.
Set 
\begin{align*}
&\overline u^\ell:=S(q')(u\times \ell)\in Z(Q'_\ell/P_\ell,\ast), 
\\
&
V_\ell:=\{[\overline u^\ell]\colon u\in U\}\subset H(Q'_\ell/P_\ell,\ast).
\end{align*}
Since the bottom arrow in the commutative diagram
\[
\begin{tikzcd}
H(Q'_\ell,P_\ell) \ar[r,"H(q')"] 
\ar[d,"H(e)"']
&
H(Q'_\ell/P_\ell,\ast)
\ar[d,"H(j)"]
\\
H(Q_\ell,P_\ell) \ar[r,"H(q)","\cong"']
&
H(Q_\ell/P_\ell,\ast)
\end{tikzcd}
\]
is an isomorphism, $[\overline u_i^\ell]$ is a basis of $V_\ell$ and 
\[
H(j)|_{V_\ell}\colon V_\ell\to H(Q_\ell/P_\ell,\ast)
\]
is an isomorphism.
Let $I_\ell$ denote the inclusion map $V_\ell\hookrightarrow H(Q'_\ell/P_\ell,\ast)$.
Define
\[
\Gamma_\ell:=H(j')I_\ell H(j)|_{V_\ell}^{-1}H(q) \colon H(Q_\ell,P_\ell)\to H(Q'_\ell/P'_\ell,\ast).
\]
For $u\in U$,
\begin{equation}
\label{eq:i_ell}
I_\ell H(j)|_{V_\ell}^{-1}H(q)[u\times \ell]=[\overline u^\ell]=[S(q')(u\times \ell)],
\end{equation}
hence
\begin{equation}
\label{eq:gamma_ell}
\Gamma_\ell [u\times \ell]=[S(j'\circ q')(u\times \ell)].
\end{equation}

\subsection{Proof of Proposition~\ref{prop:main_qp}\,(a)}
\label{ssec:proof_a}
Define 
\begin{align*}
&
F\colon (C_{[0,1]},B_{[0,1]})\to (Q'_1/P'_1,\ast),
\\
&
F(y):=\begin{cases}
\psi_{(Q',P')}(j'(q'(y)),\tau_1(y)),&\text{if $\tau_1$ is defined at $y$},
\\
\ast,&\text{otherwise},
\end{cases}
\end{align*}
hence $F(y)=j'(q'(y))$ for $y\in C_1$.
The map is continuous. Indeed, it is continuous on the set of points in which $\tau_1$
defined and that set contains 
$\clos(Q'\setminus P')_{[0,1]}$ by 
the definition of $P'$. 
On the other hand, if $y\in P'_{[0,1]}$ then
$F(y)=\ast$ independently of whether $\tau_1$ is defined or not defined at $y$,
hence the continuity of $F$ follows.
\par
Now we are able to prove the required equation
\begin{equation}
\label{eq:prop_a}
H(\Psi'_1)\Gamma_0=\Gamma_1A.
\end{equation}
Let $u\in U$. It follows by \eqref{eq:gamma_ell},
\[
H(\Psi'_1)\Gamma_0[u\times 0]=
[S(\psi_{(Q',P')}(\cdot,\tau_1(\cdot))\circ j'\circ q')(u\times 0)]
=[S(F)(u\times 0)],
\]
and
\[
\Gamma_1A[u\times 0]=[S(j'\circ q')(Au\times 1)]=[S(F)(Au\times 1)].
\] 
Since
\[
S(F)(Au\times 1)-S(F)(u\times 0)=\partial S(F)c+S(F)b\in B(Q'_1/P'_1,\ast),
\]
where the chains $c$ and $b$ are associated to the contiguous pair of cycles $(u,Au)$,
the equation \eqref{eq:prop_a} holds.\qed

\subsection{Auxiliary maps}
\label{ssec:aux}
At first we establish two facts concerning the access time maps.

\begin{lemma}
\label{lem:tau_defined}
\
\begin{itemize}
\item[(a)]
If $y\in Q_0$ and $\tau_k$ is defined at $y$ then $\kappa h<\tau_k(y)$.
\item[(b)]
If $y\in Q_k$ and $\tau_{k+m}$ is defined at $y$ then $\mu h<\tau_{k+m}(y)$.
\end{itemize}
\end{lemma}
\begin{proof}
This is an immediate consequence of 
the choice of $k$ and $m$.
\end{proof}

\begin{lemma}
\label{lem:tau_not_defined}
If for $y\in C_{[0,k]}$ the map $\tau_k$ is not defined at $y$
then there exist $t>h$,
$\epsilon >0$, and 
a neighborhood
$W$ of $y$ such that 
\begin{align}
&\label{eq:psiy}
\psi(W,[0,t+h+\epsilon])\cap Y_k=\emptyset,
\\
&\label{eq:psiq}
\psi(W,[t,t+h+\epsilon])\cap Q=\emptyset.
\end{align}
\end{lemma} 
\begin{proof}
Since $Q$ and $Y_k$ are closed, it is enough to prove the equations hold if $W$ is replaced
by the point $y$. Set
\[
t_y=\sup\{t>0\colon \psi(y,t)\in Q\},
\]
hence $t_y\geq h$.
Moreover, $t_y<\infty$, since otherwise the omega-limit set $\omega(y)$ has a nonempty
intersection with
the compact set $Q_{[0,k]}$
contradictory to Lemma~\ref{lem:no_omega}.
Since $\psi(y,t_y)\in Q$, the interval $[t_y+\epsilon,t_y+h+2\epsilon]$ is contained in
the domain of $\psi(y,\cdot)$ for a sufficiently small $\epsilon>0$ and the result is proved
for $t:=t_y+\epsilon$.
\end{proof}

For $y\in Y_{[0,k]}$ such that $\tau_k$ is defined at $y$ associate the numbers $\overline s_y\in [0,h)$ and $\overline\ell_y\in\mathbb N$ such
that
$\tau_k(y)=\overline s_y+\overline \ell_y h.$
Define a map
\begin{align*}
&
G\colon (C_{[0,k]},B_{[0,k]})\to (Q_k'/P_k,\ast),
\\
&G(y):=\begin{cases}
g^{\overline\ell_y}_{(Q',P)}(q'(\psi(y,\overline s_y))),&\text{if $\tau_k$ is defined at $y$},
\\
\ast,&\text{otherwise,}
\end{cases}
\end{align*}
thus $G(y)=q'(y)$ for $y\in C_k$. 
The map $G$ is continuous, which can be proved as follows. If $\tau_k$ is defined at $y$, it is also
defined in a neighborhood of $y$. Clearly, $G$ is continuous at $y$ if $\overline s_y>0$.
Assume $\overline s_y=0$, hence $g(y)\in Q'$. 
Let $y_n\to y$. If $\overline \ell_{y_n}=\overline \ell_y$ then $G(y_n)\to G(y)$. 
Let
$\overline s_y\to 1$. Since
 $\tau_k(y_n)$ is close to $\tau_k(y)$, $\overline\ell_{y_n}=\ell_y-1$  which implies
\[
G(y_n)=g_{(Q',P)}^{\overline \ell_y-1}(q'(\psi(y_n,\overline s_{y_n})))\to 
g^{\overline\ell-1}_{(Q',P)}(q'(g(y))) =
G(y).
\]
Assume that $\tau_k$ is not defined at $y$, hence
$G(y)=\ast$.
Let $t$, $\epsilon$, and $W$ be associated to $y$
by
Lemma~\ref{lem:tau_not_defined}. Let $y'\in W$ and $y'\neq y$. 
If $\tau_k$ is not defined at $y'$ then $G(y')=\ast$.
Finally, assume that $\tau_k$ is defined at $y'$.  
By \eqref{eq:psiy}, 
\[
t+h+\epsilon<\tau_k(y'),
\]
hence
$\overline\ell_{y'}\geq 2$ and,  
as a consequence of \eqref{eq:psiq}, there is an $\ell<\overline \ell_{y'}$ such that
at least one of the points
$g^{\ell}(\psi(y',\overline s_{y'}))$ and 
$g^{\ell+1}(\psi(y',\overline s_{y'}))$ does not belong to $Q$, hence 
$
G(y')=\ast
$
and the proof of the continuity of $G$ is finished.
\par 
Define
\[
R\colon Q'_0\to Q'_k/P_k,\quad
R(y):=
\begin{cases}
G(y),&\text{if $y\in C_0$,}
\\
\ast&\text{otherwise.}
\end{cases}
\]
The definition is correct. Indeed, since
$\clos(Q'_0\setminus C_0)\subset P'_0$ by \eqref{eq:psi_q-p}, it suffices to prove that $G(y)=\ast$ for every $y\in P'_0\cap C_0$.
This holds true if  $\tau_k$ is not defined at $y$. Assume therefore that $\tau_k$ is defined at $y$, thus
$\overline \ell_y\geq \kappa$ 
by Lemma~\ref{lem:tau_defined}\,(a). Since $\psi(y,\overline s_y)\in P'$, 
the equation
 \eqref{eq:choice_k} implies that
$g^\kappa_{(Q',P)}(q'(\psi(y,\overline s_y)))=\ast$, hence again $G(y)=\ast$.
Since $G$ is a continuous map, $R$ is also continuous. Since $R$ is equal to $\ast$ on $P'_0$, there
is a continuous map $r\colon Q'_0/P'_0\to Q'_k/P_k$ such that 
\[
r\circ j'\circ q'=R.
\]
It follows, in particular,
\begin{equation}
\label{eq:rjqu=gu}
S(r\circ j'\circ q')(u\times 0)=S(G)(u\times 0)
\end{equation}
for every $u\in U$. We assert that
\begin{equation}
\label{eq:jr=psi}
j'\circ r=\Psi'_k\colon Q'_0/P'_0\to Q'_k/P'_k.
\end{equation}
Indeed, let $y\in Q'_0\setminus P'_0$. Assume that $j'(r(y))\in Q'_k\setminus P'_k$, hence
\begin{align*}
&
g^{\overline \ell_y}(\psi(y,\overline s_y))\in Q'\setminus  P',
\\
&
g^\ell(\psi(y,\overline s_y))\in Q'\setminus P
\end{align*}
for all $0\leq \ell<\overline \ell_y$. 
This is possible only if $g^\ell(\psi(y,\overline s_y))\in Q'\setminus P'$ for
all $0\leq \ell\leq \overline\ell_y$ and the latter is equivalent to
$\psi(y,[0,\tau_k(y)])\subset Q'\setminus P'$ by the inclusion
\eqref{eq:psi_q-p} and the positive invariance of $P'$ in $Q'$ with respect
to $\psi$, hence $\Psi_k'(y)\in Q'_k\setminus P'_k$.
Since the inverse implication is trivial, the equation \eqref{eq:jr=psi} is proved.

\par
By \eqref{eq:q'q} and \eqref{eq:choice_m}, the following map is correctly defined:
\[
D\colon Q\to Q'/P,
\quad
D(y):=\begin{cases}
q'(g^\mu (y)),&\text{if $g(y),\ldots, g^\mu(y)\in Q$,}
\\
\ast,&\text{otherwise.}
\end{cases}
\]
This map is continuous.  Indeed,  let $y\in Q$.  If 
$g^\ell(y)\notin Q$ for some $1\leq \ell\leq \mu$ then
also $g^\ell(y')\notin Q$ for every $y'$ in a neighborhood of $y$,  hence
$D(y')=\ast$.  Assume that $g(y),\ldots,g^\mu(y)\in Q$ and let $y_n\to y$.
If $g(y_n),\ldots,g^\mu(y_n)\in Q$ then $D(y_n)\to D(y)$.  Assume that
$g(y_n),\ldots,g^{\ell-1}(y_n)\in Q$ and $g^\ell(y_n)\notin Q$ for some
$1\leq \ell\leq \mu$,  hence $D(y_n)=\ast$.  
Since $g^\ell (y_n)\in g(Q)\setminus Q$,  the limit
$g^\ell(y)$ of the sequence $g^\ell(y_n)$ is contained in $P$ (because $(Q,P)$ is an index pair for $g$),  
which implies $D(y)=\ast$ and the proof of
the continuity of $D$ is finished.
\par
Since $P$ is positively invariant in $Q$ with respect to $g$, 
the map $D$ factors through $Q/P$. It follows 
by Lemma~\ref{lem:tau_defined}\,(b) that the corresponding  map $Q/P\to Q'/P$ restricts to
\[
d\colon Q_k/P_k\to Q'_{[k,k+m]}/P_{[k,k+m]},
\]
i.e. $d(y)=D(y)$ if $y\in Q_k\setminus P_k$.  Moreover, again by \eqref{eq:q'q},
\begin{equation}
\label{eq:dj}
d\circ j=g^\mu_{(Q',P)}\colon Q'_k/P_k\to Q'_{[k,k+m]}/P_{[k,k+m]}.
\end{equation}
\par
Finally, define
\begin{align*}
&J\colon Q'_{[k,k+m]}/P'_{[k,k+m]}\to Q'_{k+m}/P'_{k+m},
\\
&
J(y):=
\begin{cases}
\psi_{(Q',P')}(y,\tau_{k+m}(y)),&\text{if $y\in Q'\setminus P'$ and $\tau_{k+m}$ is defined at $y$,}
\\
\ast,&\text{otherwise.}
\end{cases}
\end{align*}
The map $J$ is continuous. Indeed, it suffices to check only the case 
where $\tau_{k+m}$ is not defined at 
$y\in Q'\setminus P'$, $y_n\to y$, and  $\tau_{k+m}$ is defined at $y_n$. Let
\[
t_y:=\inf\{t>0\colon \psi(y,t)\notin \clos(Q'\setminus P')\}.
\]
There exists an $\epsilon>0$ such that 
for almost all $n$,
$t_y+\epsilon <\tau_{k+m}(y_n)$ and 
$\psi(y_n,t_y+\epsilon)\notin\clos(Q'\setminus P')$, hence $J(y_n)=\ast$ and the continuity of $J$
is proved. By definition,
\begin{equation}
\label{eq:j=psi}
J|_{Q'_k/P'_k}=\Psi'_m.
\end{equation}
 Moreover, 
if $y\in Q'_{[k,k_m]}\setminus P'_{[k,k+m]}$ and $\tau_{k+m}$ is
 defined at $y$ then for $t\leq \tau_{k+m}(y)$,
\[ 
 J(y)=J(\psi_{(Q',P')}(y,t)),
 \]
hence, by Lemma~\ref{lem:tau_defined}\,(b), 
\begin{equation}
\label{eq:j=jg}
J|_{Q'_k/P'_k}=J\circ g^\mu_{(Q',P')}|_{Q'_k/P'_k}.
\end{equation}

\subsection{Proof of Proposition~\ref{prop:main_qp}\,(b)}
\label{ssec:proof_b}
We will prove that
\[
\Delta:=A^mH(q)^{-1}H(j)H(r)\colon H(Q'_0/P'_0,\ast)\to H(Q_{k+m},P_{k+m})
\]
is the required map, i.e.
\begin{align}
\label{eq:delta_up}
&
\Delta\Gamma_0=A^{k+m},
\\
\label{eq:delta_down}
&
\Gamma_{k+m}\Delta=H(\Psi'_{k+m}).
\end{align}
In order to prove \eqref{eq:delta_up} we will show that the diagram
\begin{equation}
\label{eq:diagram}
\begin{tikzcd}[column sep=huge,row sep=large]
H(Q_0,P_0) \ar[r,"A^k"]
\ar[d,"\Gamma_0"] 
&
H(Q_k,P_k) \ar[d,"I_kH(j)|_{V_k}^{-1}H(q)"]
\\
H(Q'_0/P'_0,\ast) \ar[r,"H(r)"]
&
H(Q'_k/P_k,\ast) 
\end{tikzcd}
\end{equation}
commutes. To this purpose we should prove that
\begin{equation}
\label{eq:rg=l}
H(r)\Gamma_0=\Lambda,
\end{equation}
where 
$\Lambda$ denotes the composition of the top and the right-hand side arrows, i.e.
\[
\Lambda:=I_k H(j)|_{V_k}^{-1}H(q)A^k\colon H(Q_0,P_0)\to H(Q'_k/P_k,\ast).
\]
Let $u\in U$. By \eqref{eq:gamma_ell}
and \eqref{eq:rjqu=gu},
\[
H(r)\Gamma_0[u\times 0]=[S(r\circ j'\circ q')(u\times 0)]=[S(G)(u\times 0)].
\]
On the other hand, by \eqref{eq:i_ell},
\[
\Lambda[u\times 0]=[S(q')(A^ku\times k)]=[S(G)(A^ku\times k)].
\]
Since 
\[
S(G)(A^ku\times k)-S(G)(u\times 0)=S(G)c+S(G)b\in B(Q'_k/P_k,\ast)
\]
for chains $c$ and $b$ associated to the contiguous pair of cycles $(u,A^ku)$,
the left and the right hand sides of \eqref{eq:rg=l} coincide and thus the diagram \eqref{eq:diagram} commutes.
The commutativity of \eqref{eq:diagram} immediately implies
the equation \eqref{eq:delta_up}. Indeed,
\[
\Delta\Gamma_0=A^m H(q)^{-1}H(j)I_k H(j)|_{V_k}^{-1}H(q)A^k=A^{k+m}.
\]
It remains to prove the equation \eqref{eq:delta_down}. To this end 
let us take an arbitrary homology class $[z]\in H(Q'_k/P_k,\ast)$, where $z$ is cycle
in $Z(Q'_k/P_k,\ast)$.
Since 
\[
H(q)\colon H(Q_k,P_k)\to (Q_k/P_k,\ast)
\]
is an isomorphism,
there exists a $u\in U$ such that
\begin{equation}
\label{eq:hjz=hdu}
H(j)[z]=H(q)[u\times k],
\end{equation}
hence $H(j)[z]=H(j)[\overline u^k]$ and therefore there exist
chains $\eta\in S(Q_k/P_k)$ and $\xi\in S(\ast)$ satisfying
\begin{equation}
\label{eq:sjz+sju}
S(j)z-S(j)\overline u^k=\partial \eta +\xi.
\end{equation}
Since the diagram
\[
\begin{tikzcd}[column sep=large]
Q'/P \ar[r,"g_{(Q',P)}"] 
\ar[d,"j'"]
&
Q'/P \ar[d,"j'"]
\\
Q'/P' \ar[r,"g_{(Q',P')}"]
&
Q'/P'
\end{tikzcd}
\]
commutes 
(because  $P'$ is positively invariant in $Q'$)
and the equation \eqref{eq:dj} holds, applying $S(j'\circ d)$ to the equation \eqref{eq:sjz+sju} one gets
\[
S(g^\mu_{(Q',P')})S(j')z-S(g^\mu_{(Q',P')})S(j')\overline u^k=\partial S(j'\circ d)\eta+\xi,
\]
hence, by \eqref{eq:j=jg},
\[
S(J)S(j')z-S(J)S(j')\overline u^k=\partial S(J\circ j'\circ d)\eta+\xi
\]
and therefore the equation \eqref{eq:j=psi} implies
\begin{equation}
\label{eq:hz=hu}
H(\Psi'_m)H(j')[z]=H(\Psi'_m)H(j')[\overline u^k].
\end{equation}
By definition, 
\[
H(j')[\overline u^k]=\Gamma_k[u\times k],
\] 
and by an iterated application of \eqref{eq:prop_a},
\[
H(\Psi'_m)\Gamma_k[u\times k]=\Gamma_{k+m}A^m[u\times k],
\]
hence the equations \eqref{eq:hjz=hdu} and \eqref{eq:hz=hu} imply
\[
H(\Psi'_m)H(j')[z]=\Gamma_{k+m}A^mH(q)^{-1}H(j)[z].
\]
The latter equation together with the equation \eqref{eq:jr=psi} provide the commutativity of the diagram
\[
\begin{tikzcd}[column sep=2cm]
&
H(Q_k,P_k) \ar[r,"A^m"]
&
H(Q_{k+m},P_{k+m}) \ar[dd,"\Gamma_{k+m}"]
\\
&
H(Q'_k/P_k,\ast) \ar[u,"H(q)^{-1}H(j)"]
\ar[d,"H(j')"']
&
\\
H(Q'_0/P'_0,\ast) 
\ar[r,"H(\Psi'_k)"]
\ar[ru,"H(r)"]
&
H(Q'_k/P'_k,\ast) \ar[r,"H(\Psi'_m)"]
&
H(Q'_{k+m},P'_{k+m})
\end{tikzcd}
\]
Thus
\[
\Gamma_{k+m}\Delta=H(\Psi'_m)H(\Psi'_k)
\]
and therefore the equation \eqref{eq:delta_down} holds. This concludes
the proof of Proposition~\ref{prop:main_qp} as well as the proof
of  Theorem~\ref{thm:main}. \qed

\appendix
\section{}
\label{sec:app}

\subsection{Proof of Proposition~\ref{prop:leray}}
\label{ssec:a1}
Since $\operatorname{gker}(\alpha)=\operatorname{gker}(\alpha^r)$, $\overline{\alpha^r}=\overline\alpha^r$ and 
since $\operatorname{gim}(\overline\alpha)=\operatorname{gim}(\overline \alpha^r)$, 
\[
R(\alpha^r)=(R\alpha)^r,
\quad
RV_{\alpha^r}=RV_\alpha
\]
and the same equations hold with $\alpha$ replaced by $\beta$.
Moreover, the maps 
\begin{align*}
&
\widehat\gamma\colon RV_\alpha\to RV_\beta,\quad 
\widehat\gamma(v+\operatorname{gker}(\alpha)):=\gamma(v)+\operatorname{gker}(\beta),
\\
&
\widehat\delta\colon RV_\beta\to RV_\alpha,\quad
\widehat\delta(v+\operatorname{gker}(\beta)):=\delta(v)+\operatorname{gker}(\alpha)
\end{align*}
are well defined and the diagrams
\[
\begin{tikzcd}[column sep=large]
RV_\alpha \ar[r,"R\alpha"] \ar[d,"\widehat\gamma"]
&
RV_\alpha \ar[d,"\widehat\gamma"]
\\
RV_\beta \ar[r,"R\beta"]
&
RV_\beta
\end{tikzcd}
\qquad
\begin{tikzcd}[column sep=large]
RV_\alpha \ar[r,"(R\alpha)^r"] \ar[d,"\widehat\gamma"]
&
RV_\alpha \ar[d,"\widehat\gamma"]
\\
RV_\beta \ar[r,"(R\beta)^r"] \ar[ru,"\widehat\delta"]
&
RV_\beta
\end{tikzcd}
\]
commute.
Since $(R\alpha)^r=\widehat\delta\widehat\gamma$ is an automorphism,
$\widehat\gamma$ is a monomorphism,
and since 
$(R\beta)^r=\widehat\gamma\widehat\delta$ is an automorphism,
$\widehat\gamma$ is an epimorphism.
Therefore $\widehat\gamma$ is an isomorphism and the result follows. 
\qed

\subsection{Proof of Proposition~\ref{prop:s_invariant}}
\label{ssec:s_inv}
As it was already mentioned, 
the equivalence of (c) and (d) is a general result, valid for any  
local dynamical system $\Phi$ on a locally compact space $X$. Here is a proof:
\par\noindent
(c) $\Rightarrow$ (d). We can assume $t>0$. Since $S$ is invariant with respect to $\Phi$,
it is also invariant with respect to $\Phi_t:=\Phi(\cdot,t)$. If $N$ is an isolating neighborhood
of $S$ in the system $\Phi$ then there is a compact neighborhood
$W$ of $S$ such that $\Phi(W,[-t,t])\subset N$, hence $W$ is an isolating neighborhood of $S$ in $\Phi_t$.
\par\noindent
(d) $\Rightarrow$ (c). Let $t\neq 0$, let $S$ be an isolated invariant set with respect
to $\Phi_t$, and let $N$ be its isolating neighborhood. There is an $\epsilon>0$
such that $\Phi(S,[-\epsilon,\epsilon])\subset N$. Since $\Phi_\epsilon(S)$ is invariant with respect to
$\Phi_t$, $\Phi_\epsilon (S)\subset S$, hence $\Phi(S,[-\epsilon,2\epsilon])\subset N$, etc. Consequently,
 $\Phi(S,\mathbb R)=S$,
i.e. $S$ is invariant for $\Phi$. Since each invariant set with respect to $\Phi$
contained in $N$ is also invariant for $\Phi_t$, $N$ is also an isolating
neighborhood of $S$ with respect to $\Phi$.
\par 
Now we assume that $\Phi$ is the translation system associated to $\phi$. In order
to finish the proof it suffices to prove that (d) is equivalent to (e). We follow the notation 
from Subsection~\ref{ssec:poincare}.
\par\noindent
(d) $\Rightarrow$ (e). Let $N$ be an isolating neighborhood of $S$ with respect to $\Phi_1$. 
Then $N_a$ is an isolating neighborhood of $S_a$ for $\Phi_1$ restricted to $X_a\times X_a$. Let $t\in [0,1)$,
hence $\Phi_t(S_a)\subset S_{a+t}$ and $\Phi_{-t}(S_{a+t})\subset S_a$ by the the invariance
of $S$ with respect to $\Phi$ stated in the implication (d) $\Rightarrow$ (c).
Thus $\Phi_t(S_a)=S_{a+t}$.
\par\noindent
(e) $\Rightarrow$ (d). Let $a,t\in \mathbb R$. Since the maps $\Phi_t$ and $\Phi_{-t}$ are mutually inverse
homeomorphisms in neighborhoods of $S_a$ and $S_{a+t}$, respectively, $S_a$ is isolated invariant
if and only $S_{a+t}$ is isolated invariant in $\Phi_1$.  For $t\in [0,1)$ let $N_t$ be an isolating neighborhood
of $S_t$. There exists an $\epsilon_t>0$ such that $\zeta(N_t\times (t-\epsilon_t,t+\epsilon_t))_{t'}$ 
is an isolating neighborhood of $S_{t'}$ if $|t'-t|<\epsilon_t$. Let $t_1,\ldots,t_k\in (0,1)$ be such that
\[
[0,1]\subset
(-\epsilon_0,\epsilon_0)\cup (1-\epsilon_0,1+\epsilon_0)+\bigcup_{i=1,\ldots,k} 
(t_i-\epsilon_{t_i},t_i+\epsilon_{t_i}).
\]
Set $t_0:=0$. Then the union of the sets $\zeta(N_{t_i}\times (t_i-\epsilon_{t_i},t_i+\epsilon_{t_i}))$ for
$i=0,\ldots,k$ is an isolating neighborhood of $S$ with respect to $\Phi_1$, hence (d) is satisfied for $t=1$.
Thus (d) is satisfied for an arbitrary $t$ by
the equivalence of (c) and (d).\qed

\subsection{Proof of Proposition~\ref{prop:index_pairs}}
\label{ssec:a3}
At first we prove that
the sets $K$, $M$, and $N'$ are compact.
We begin with a proof of the compactness of $K$.
Since $L$ is positively invariant in $N$ with respect to $\phi^h$,
\[
K=\bigcup\{\sigma\subset N\colon \text{$\sigma$ is a segment with length in $[0,h]$ starting at a point in $L$}\}.
\]
Let $x_n\in K$, $n\in\mathbb N$. There are $y_n\in L$ and $t_n\in [0,h]$ such that $x_n=\phi(y_n,t_n)$
and $\phi(y_n,[0,t_n])\subset N$. We can assume $y_n\to y\in L$ and $t_n\to t\in [0,h]$ as $n\to \infty$;
it follows
$\phi(y,[0,t])\subset N$. Therefore $x_n\to \phi(y,t)\in K$ as $n\to\infty$
and the compactness of $K$ follows.
Clearly, an analogous argument implies that the set $M$
is compact, hence $N'$ is a compact set as well.

\begin{lemma}
\label{lem:step_3}
If $x\in M$, $\phi^h(x)\in N$, and $\phi(x,[0,h])\not\subset N$ then $\phi^h(x)\in K$.
\end{lemma}
\begin{proof}
By assumption 
 there is a number $s$ such that
\begin{equation}
\label{eq:s}
s\in (0,h],\quad
\phi(x,[s,h])\subset N,\quad 
\exists \epsilon_n\in (0,s),\ \epsilon_n\to 0\colon
\phi(x,s-\epsilon_n)\notin N
\end{equation}
Moreover,
there are $y\in N$ and $t\in (0,h]$ such that $\phi(y,[0,h])\subset N$ and 
$x=\phi(y,t)$, hence
$\phi(x,[0,h-t])\subset N$.
Then
\[
h<s+t-\epsilon_n<s+t\leq 2h
\]
for every $n$,
since otherwise
\[
\phi(x,[0,h-t])\cup \phi(x,[s-\epsilon_n,h])=\phi(x,[0,h])\subset N;
\]
a contradiction. It follows 
\[
z_n:=\phi(y,s+t-h-\epsilon_n)\in N,\quad z_n\to z:=\phi(y,s+t-h), 
\]
and, by \eqref{eq:s},
\[
\phi^h(z_n)=\phi(x,s-\epsilon_n)\notin N,
\quad
\phi^h(z)=\phi(x,s)\in N.
\]
Since 
$L$ is an exit set of $N$, $\phi(x,s)\in L$ and therefore $\phi^h(x)\in K$.
\end{proof}

\begin{lemma}
\label{lem:step_4}
If $x\in K$ and $\phi^h(x)\in N$ then $\phi^h(x)\in K$.
\end{lemma}
\begin{proof}
Let $y\in L$, $x=\phi(y,t)$, and $\phi(y,[0,t])\subset N$. Since $L$ is positively invariant in $N$,
we can assume $t\in (0,h)$. If $\phi(x,[0,h])\subset N$ then also $\phi(y,[0,h+t])\subset N$, hence
$\phi^h(x)=\phi(y,t+h)\in K$. Assume that $\phi(x,[0,h])\not\subset N$. Let $s$ satisfy \eqref{eq:s}. 
If $h<s+t$ then also $h<s+t-\epsilon_n$ for sufficiently large $n$, hence
\[
\phi(y,s+t-h-\epsilon_n)\in \phi(y,(0,t))\subset N
\]
and, by the argument in the proof of Lemma~\ref{lem:step_3}, we conclude that $\phi^h(x)\in K$. Finally, assume that $s+t\leq h$.
Then
\[
\phi^h(y)=\phi(x,h-t)\in \phi(x,[s,h))\subset N,
\]
hence $\phi^h(y)\in L$ because $L$ is positively invariant in $N$. Therefore 
\[
\phi^h(x)\in \phi(\phi^h(y),[0,t])\subset K
\]
and the proof is finished.
\end{proof}

\begin{proof}[Proof of {\rm (b)}]
Let $x\in N'$. Assume first that  
$x\in M$. If 
$
\phi(x,[0,h])\subset N
$
then $\phi^h(x)\in M\subset N'$. If, on the other
hand, $\phi(x,[0,h])\not\subset N$ then $\phi^h(x)\in K\subset N'$ by Lemma~\ref{lem:step_3}. Now let $x\in K$.
Then $\phi^h(x)\in K\subset N'$ by Lemma~\ref{lem:step_4}.
\end{proof}

\begin{proof}[Proof of {\rm (a)}]
Let $W$ be a neighborhood of $S$ such that $\phi(W,[0,h])\subset \interior(N\setminus L)$.
It follows $S\subset \interior N'$, hence also
$S\subset \interior(N'\setminus L)$ and therefore
\[
S\subset \Inv\clos(N'\setminus L)\subset \Inv\clos(N\setminus L)=S
\]
which means the pair $(N',L)$ is isolating. 
Since $N'\subset N$, $L$ is positively invariant in $N'$. Finally, let $x_n\in N'$ be such that
$\phi^h(x_n)\notin N'$. Assume that
 $x_n\to x$ and $\phi^h(x)\in N'$. Since $\phi^h(x_n)\notin N$ by the conclusion (b) 
 and $L$ is an exit set of $N$,
 $\phi^h(x)\in L$ and therefore $L$ is an exit set of $N'$ as well.
\end{proof}

\begin{proof}[Proof of {\rm (c)}]
By Proposition~\ref{prop:s_invariant}, $S$ is an isolated invariant set with respect to $\phi$.
Since $S\cap K=\emptyset$, $K$ is closed, and the conclusion (a) holds,
\[
S\subset \interior(N'\setminus L)\cap (X\setminus K)= \interior(N'\setminus K)
\]
and therefore
\[
\Inv_\phi \clos(N'\setminus K)\subset \Inv_{\phi^h}(\clos(N'\setminus L)=S\subset \interior(N'\setminus K),
\]
i.e. the pair $(N',K)$ is isolating. It follows directly by
definition that $K$ is positively invariant in $N'$ and it remains to prove
that $K$ is an exit set of $N'$. To this reason assume
that $t\geq 0$ and
\[
\phi(x,[0,t])\subset N',\quad \phi(x,t+\epsilon_n)\notin N'
\]
for a sequence $\epsilon_n$ such that $\epsilon_n>0$, $\epsilon_n\to 0$.
If $x\in K$ then obviously $\phi(x,t)\in K$. In the other
case $x\in M$, hence $x=\phi(y,h-t)$ for some $y\in N$ and
$\phi(y,[0,h])\subset N$.
Set $y_n:=\phi(y,\epsilon_n)$. Clearly, $y_n\in N'$ for $n$ large enough
and $\phi^n(y_n)\notin N'$,
hence
$\phi^h(y_n)\notin N$ by the conclusion (b). 
Since $L$ is an exit set of $N$ with respect
to $\phi^h$, $\phi(x,t)=\phi^h(y)\in L$ and the required property 
of the pair $(N,L)$ is proved.
\end{proof}

\subsection{Proof of Proposition~\ref{prop:a4}}
\label{ssec:a4}
 As one can check directly, $(N',L')$ is an index pair for $S$ with respect to  
$\phi$
(see also \cite[Rem.\,5.4]{salamon}).
By Proposition~\ref{prop:index_pairs}\,(a), for the system $\phi^h$ the set  
$S$ equal to the invariant part of $\clos(N'\setminus L')$ and $L'$ is an
exit set of $N'$, hence in order to prove that $(N',L')$ is an index
pair for $S$ with respect to $\phi^h$ it remains to prove that 
$L'$ is positively invariant in $N'$. Let $x\in L'$ and assume that $\phi^h(x)\in N'$. 
If $\phi(x,[0,h])\subset N'$ then clearly $\phi^h(x)\in L'$.
If, in the other hand, $\phi(x,[0,h]\not\subset N'$ then Lemmas~\ref{lem:step_3} and \ref{lem:step_4}
guarantee that $\phi^h(x)\in K\subset L'$ and the result follows.\qed

\end{document}